\documentstyle[12pt]{article}

\def\be{\begin{equation}}
\def\ee{\end{equation}}
\def\uom{\underline{\omega}}
\def\IP{\hbox{\rm I\kern -1.6pt{\rm P}}}
\def\IC{{\hbox{\rm C\kern-.58em{\raise.53ex\hbox{$\scriptscriptstyle|$}}
    \kern-.55em{\raise.53ex\hbox{$\scriptscriptstyle|$}} }}}
\def\IN{\hbox{I\kern-.2em\hbox{N}}}
\def\IR{\hbox{\rm I\kern-.2em\hbox{\rm R}}}
\def\ZZ{\hbox{{\rm Z}\kern-.33em{\rm Z}}}
\def\IT{\hbox{\rm T\kern-.38em{\raise.415ex\hbox{$\scriptstyle|$}} }}
\def\ang{\hbox{$<$\kern-.42em\hbox{\rm )} }}
\def\ggm{{G}/\Gamma}

\begin{document}

\title{Dynamical Borel-Cantelli lemmas for Gibbs measures}
\author{N. Chernov$^\ast$ and D. Kleinbock$^\dagger$
}
\date {\it To appear in Israel J.~Math.}
\maketitle

\begin{abstract}
Let $T:\, X\mapsto X$ be a deterministic dynamical system preserving a
probability measure $\mu$. A dynamical Borel-Cantelli lemma
asserts that for certain sequences of subsets $A_n\subset X$ and
$\mu$-almost every point $x\in X$ the inclusion $T^nx\in A_n$
holds for infinitely many $n$. We discuss here systems which are
either symbolic (topological) Markov chain or Anosov
diffeomorphisms preserving Gibbs measures. We find sufficient
conditions on sequences of cylinders and rectangles, respectively,
that ensure the dynamical Borel-Cantelli lemma.
\end{abstract}

\newtheorem{theorem}{Theorem}[section]
\newtheorem{lemma}[theorem]{Lemma}
\newtheorem{sublemma}[theorem]{Sublemma}
\newtheorem{proposition}[theorem]{Proposition}
\newtheorem{corollary}[theorem]{Corollary}

\renewcommand{\theequation}{\arabic{section}.\arabic{equation}}

\footnotetext{$^\ast$ Partially supported by NSF grant
DMS-9732728.}

\footnotetext{$^\dagger$ Partially supported by NSF grant
DMS-9704489.}

\section{Introduction}
\label{secI}
\setcounter{equation}{0}

Let $T: X\mapsto X$ be a transformation preserving a probability
measure $\mu$. We use notation $\mu(f):=\int f\, d\mu$ for
integrable functions $f$ on $X$.

Let $A_n\subset X$ be a sequence of measurable sets. Put
$B_n=T^{-n}A_n$ and consider the set $$ \limsup_n
B_n:=\cap_{m=1}^{\infty}\cup_{n=m}^{\infty}B_n$$ of points which belong
to infinitely many $B_n$. A classical
Borel-Cantelli lemma in probability theory states:

\begin{lemma}[Borel-Cantelli]
{\rm (i)} If $\sum\mu(B_n)<\infty$, then $\mu(\limsup_n
B_n)=0$, i.e.~almost
every point $x\in X$ belongs to finitely many $B_n$.\\ {\rm (ii)}
If $\sum\mu(B_n)=\infty$ and $B_n$ are independent, then
$\mu(\limsup_n
B_n)=1$, i.e.~almost every point $x\in X$ belongs to infinitely
many $B_n$.
\end{lemma}

In terms of the transformation $T$, the lemma can be restated as
follows.

\begin{lemma}
{\rm (i)} If $\sum\mu(A_n)<\infty$, then for almost every point
$x\in X$ there are only finitely many $n$ such that $T^nx\in
A_n$.\\ {\rm (ii)} If $\sum\mu(A_n)=\infty$ and $T^{-n}A_n$ are
independent, then for almost every point $x\in X$ there are
infinitely many $n$ such that $T^nx\in A_n$. \label{lmBCT}
\end{lemma}

The second part of the lemma has a limited value for deterministic
dynamical systems, since one rarely works with purely independent
sets.
This paper is devoted to extensions of the second part of the
lemma to certain dynamical systems -- Anosov diffeomorphisms and
topological Markov chains.

Below we always assume that $\sum_n\mu(A_n)=\infty$.
\medskip

\noindent {\bf Definition}. A sequence of subsets $A_n\subset X$
is called a {\sl Borel-Cantelli} (BC) sequence if for $\mu$-a.e.~$x\in
X$ there are infinitely many $n$ such that $T^nx\in A_n$.
\medskip

Let
$$
    \chi_n(x):=\chi_{T^{-n}A_n}(x)
$$
be the indicator of the set $B_n=T^{-n}A_n $. We set
$$
    S_N(x):=\sum_{n=1}^N\chi_n(x)
$$
and
$$
    E_N':=\mu(S_N)=\sum_{n=1}^N\mu(A_n)\,.
$$

\noindent {\bf Definition}. A sequence of subsets $A_n\subset X$
is said to be a {\sl strongly Borel-Cantelli} (sBC) sequence if
for $\mu$-a.e.~$x\in X$ we have $S_N(x)/E_N\to 1$ as $N\to\infty$.
\medskip

A stronger version of the classical Borel-Cantelli lemma is known,
see Theorem~6.6 in \cite{Durr}:

\begin{lemma} If $\sum\mu(B_n)=\infty$ and the events $B_n$ are
independent, then $S_N(x)/E_N\to 1$ almost surely as
$N\to\infty$. Moreover, the independence requirement can be
relaxed to the pairwise independence, i.e.~it is enough to require
$\mu(B_m\cap B_n)=\mu(B_m)\mu(B_n)$ for $m\neq n$.
\end{lemma}

In particular, if $B_n=T^{-n}A_n$ are pairwise independent, then
the sequence $\{A_n\}$ is an sBC sequence.

Consider the quantity
$$
      R_{mn}:=\mu(B_m\cap B_n)-\mu(B_m)\mu(B_n)=\mu(T^{-m}A_m\cap
T^{-n}A_n)-\mu(A_m)\mu(A_n)
$$
which characterizes the dependence of $B_m$ and $B_n$.

A sufficient condition for $\{A_n\}$ to be an sBC sequence, in
terms of $R_{mn}$, was first found by W.~Schmidt, see a proof by
Sprind\v{z}uk \cite{Sp}, in the context of Diophantine
approximations. It was recently adapted to dynamical systems by
D.~Kleinbock and G.~Margulis \cite{KM}:

\noindent{\bf (SP)} Assume that
$$
    \exists C>0:\ \ \ \ \sum_{m,n=M}^NR_{mn}\leq C\cdot \sum_{n=M}^N\mu(A_n)
$$
for all $N\geq M\geq 1$.

\begin{theorem}[\cite{Sp}{\rm , Chapter I,  Lemma 10, or} \cite{KM}{\rm , Lemma 2.6}]
If the sequence $\{A_n\}$ satisfies {\rm (SP)}, then it is an sBC
sequence; moreover,  for a.e.~$x\in X$ one has
\be
     S_N=E_N+O\Big(E_N^{1/2}\log^{3/2+\varepsilon}E_N\Big )\,.
          \label{log}
\ee
\end{theorem}

W.~Philipp was first to derive the asymptotics (\ref{log}) in the context
of dynamical system, and he called it a quantitative Borel-Cantelli lemma
\cite{Ph}.\medskip

Note that there exist remarkable characterizations of
some ergodic properties of dynamical systems in terms of BC and sBC sequences.
We summarize these in the following

\begin{proposition}  Let $T$ be  a measure preserving transformation
of a probability space $(X,\mu)$. Then:\\
{\rm (i)} $T$ is ergodic $\iff$
every constant sequence $A_n\equiv A$, $\mu(A)>0$, is BC $\iff$ every
such sequence is sBC, i.e.~$S_N/E_N\to 1$ $\mu$-almost everywhere;\\
{\rm (ii)}
$T$ is weakly mixing $\iff$ every sequence $\{ A_n\}$ that only contains finitely
many distinct sets, none of them of measure zero, is BC  $\iff$ for
every such sequence one has
$S_N/E_N\to 1$ in the $L^2$ metric, i.e.~$\mu(S_N/E_N-1)^2\to 0$;\\
{\rm (iii)} $T$
is lightly mixing\footnote{ $T$ is said to be
lightly mixing (see \cite{FT}) if for every two sets $A,B$ of positive
measure one has \linebreak
$\mu(T^{-n}A \cap B) > 0$ for large enough $n$; this condition lies
strictly between mixing and weak mixing.} $\iff$ every sequence  that
only contains finitely
many distinct sets, possibly of measure zero, is BC.
\label{lmWM}
\end{proposition}

See Section~\ref{secPTMC} for the proof. Note that in part (ii),
the first equivalence  was proved by Y.~Guivarc'h and A.~Raugi
(private communication); our proof is slightly different. Part
(iii) was pointed out to us by A.~del Junco. \medskip

Note also that there exist no measure-preserving system such that
every sequence $\{ A_n\}$ that only contains two distinct sets,
one of positive measure and the other of measure zero, is sBC.
This follows from a result of U. Krengel \cite{Kr}. On the other
hand, if  $\mu$ has $K$ property, then any sequence  that only
contains finitely many sets, none of them of measure zero, is sBC
(J.-P.~Conze, private communication).


\medskip


It is important to mention that for any
(nontrivial) measure-preserving system $(X,\mu,T)$  there are
sequences of subsets
of $X$ (with divergent sum of measures) which are not  BC. More precisely,
the following is true:

\begin{proposition}  Let  $(X,\mu)$ be a probability space. If $\mu$ is nontrivial
(that is, there are sets with measure strictly between $0$ and $1$), then
for any $\mu$-preserving transformation $T$ of $X$
there exists a sequence
$\{A_n\}$ of
measurable subsets of $X$  with $\sum_{n = 1}^\infty
\mu({A}_n)= \infty$ which is not BC. Furthermore, if $\mu$ is non-atomic, then
for any $\mu$-preserving transformation $T$ of $X$
there exists a sequence
$\{A_n\}$ of
measurable subsets of $X$  with $\sum_{n = 1}^\infty
\mu({A}_n)= \infty$ such that  for
a.e.~$x\in X$  there are
 at most finitely many $n$ for which $T^n x\in A_n$.
\label{lmNBC}
\end{proposition}

See the end of  Section~\ref{secPTMC} for the proof.
With a little extra work, one can always find a non BC sequence of
sets that are nested: $A_1\supset A_2\supset\cdots$. We omit the
proof.

Observe that a non-BC sequence can be easily constructed when  $T$
is invertible:  one can simply take $A_n=T^nA$, where $0 < \mu(A)
< 1$.
 Therefore to prove the BC or sBC property for certain classes of
sequences it is necessary to impose certain restrictions on the
sets $A_n$, which, roughly speaking, guarantee that the sets $B_m$
and $B_n$ become  nearly independent for large $|m-n|$.

\medskip

The first Borel-Cantelli lemma for deterministic
dynamical systems was proved in 1969 by W.~Philipp:

\begin{theorem}[\cite{Ph}]
Assume that $T(x)=\beta x$ (mod 1) with $\beta>1$, or
$T(x)=\{1/x\}$ (the Gauss transformation) and $\mu$ is the unique
$T$-invariant smooth measure on $[0,1]$. Then any sequence
$\{A_n\}$ of subintervals (with divergent sum of measures) is an sBC
sequence, and {\rm (\ref{log})} holds.
\label{tmP}
\end{theorem}

In particular,  one can take any $x_0\in (0,1)$ and consider what
could be called
``a target shrinking
to $x_0$'' (terminology borrowed from \cite{HV}), i.e.~a sequence  of
intervals $A_n = (x_0 - r_n, x_0 +
r_n)$ with $r_n\to 0$. Then almost all orbits $\{T^nx\}$ get into infinitely many such
intervals whenever $r_n$ decays  slowly enough. This can be
thought of as a quantitative strengthening of density of almost all orbits
(cf.~the paper \cite{Bos} for a similar approach to the
rate of recurrence).

More generally, if  $X$ is a metric space (e.g.~a Riemannian
manifold), one can try to prove that
any sequence
$\{A_n\}$ of balls in $X$ is BC or  sBC; as in the example above, this
would imply that all points $x_0\in X$ can be ``well approximated'' by
orbit points $T^nx$ for almost all $x$.
D. Dolgopyat recently proved the following:

\begin{theorem}[\cite{Do}]
Let $T:X\mapsto X$ be an Anosov diffeomorphism with a smooth invariant
probability measure $\mu$. Then any sequence of round balls
(with divergent sum of measures) is sBC.
\label{tmD}
\end{theorem}

Another example of a dynamical Borel-Cantelli lemma is given in the
paper \cite{KM}, where the following theorem was essentially proved:

\begin{theorem}[\cite{KM}]
 Let $G$ be a connected semisimple center-free Lie
group  without compact factors,  $\Gamma$ an irreducible
lattice in $G$, $\mu$ the normalized Haar measure on $\ggm$, $g$ a
partially hyperbolic element of $G$, and let $T$ be the left shift
$T(x) = gx$, $x\in \ggm$. Let $\{A_n\}$ be a sequence of subsets of
$\ggm$ with divergent sum of measures and ``uniformly regular
boundaries'', namely, such that for some $\delta>
0$ and $0 < c < 1$ one has
\be
\mu(\delta{\rm -neighborhood \ of\  }\partial A_n) \le c\mu(A_n) {\rm
\ \ \ for\ 
all\  }n\,.
\label{ur}
\ee
Then there exist positive $C_1, C_2$ such that for $\mu$-a.e.~$x\in
\ggm$ one has
$$
 C_1 \le \liminf_{N\to\infty}S_N(x)/E_N \le
\limsup_{N\to\infty}S_N(x)/E_N \le C_2\,;
$$
in particular, $\{A_n\}$ is a BC sequence.
\label{tmKM}
\end{theorem}

It is shown in \cite{KM} that the above condition (\ref{ur}) is
satisfied if $\ggm$ is not compact and the sets $\{A_n\}$ are
{\em complements} of balls centered in a fixed point $x_0\in\ggm$.
This way one gets a description of {\em growth} of almost all orbits
$T^nx$ as follows: if a sequence $R_n$ increases slowly enough, then
for almost all $x$ one has dist$(x_0,T^nx) \ge R_n$ for infinitely
many $n$. This has important applications to geometry and
number theory.

When this paper was under preparation, we learned that J.-P.~Conze
and A.~Raugi \cite{CR} proved a dynamical Borel-Cantelli lemma for
certain Markov processes and one-sided topological Markov chains
with Gibbs measures.

\section{Statement of results}
\label{secSR}
\setcounter{equation}{0}

Our paper deals with Anosov diffeomorphisms and the corresponding
symbolic systems -- topological Markov chains.

Let $T:X\mapsto X$ be a transitive Anosov diffeomorphism. Let ${\cal
R}=\{R_1,\ldots,R_M\}$ be a finite Markov partition of $X$, and
$\bf A$ the corresponding transition matrix of zeroes and ones.
For definitions and basic facts on Markov partitions, see
\cite{Bo75,ChS}.

The matrix $\bf A$ is transitive, i.e.~${\bf A}^K$ is completely
positive for some $K\geq 1$. Let $\Sigma=\Sigma_{\bf A}$ be the
{\sl topological Markov chain} for $\bf A$, i.e.~a set of doubly
infinite sequences $\uom=\{\omega_i\}_{i=-\infty}^{\infty}\in
\{1,\ldots,M\}^{\ZZ}$ defined by
$$
   \Sigma=\{\uom\in\{1,\ldots,M\}^{\ZZ}:\
    {\bf A}_{\omega_i\omega_{i+1}}=1\ \ \ \forall i\in\ZZ\}\,.
$$
The set $\Sigma$ equipped with the product topology is a
compact space, and there is a left shift homeomorphism
$\sigma:\Sigma\mapsto\Sigma$ defined by
$(\sigma\uom)_i=\omega_{i+1}$. Let $\pi:\Sigma\mapsto X$ be
the projection defined by
$$
    \pi(\uom)=\cap_{i=-\infty}^{\infty}T^{-i}R_{\omega_i}\,.
$$
Then $\pi$ is a continuous surjection and
$\pi\circ\sigma=T\circ\pi$. Fix an $a\in (0,1)$ and let $d_a$ be a
metric on $\Sigma$ defined by $d_a(\uom,\uom')=a^n$ where
$n=\max\{n:\, \omega_i=\omega_i',\ \forall |i|<n\}$. It is
consistent with the product topology. The projection $\pi$ is now
H\"older continuous.

There are classes of {\em Gibbs measures} on both $X$ and $\Sigma$
defined by potential functions. For any H\"older continuous
function $\psi:\,\Sigma\mapsto\IR$ there is a unique
$\sigma$-invariant Gibbs measure $\mu_{\psi}$ on $\Sigma$. For any
H\"older continuous function $\varphi:\, X\mapsto\IR$ there is a
unique $T$-invariant Gibbs measure $\mu_{\varphi}$ on $X$. In the
latter case, the function $\psi=\varphi\circ\pi$ is H\"older
continuous on $\Sigma$, and the measure $\mu_{\psi}$ projects to
$\mu_{\varphi}$ in the sense that $\pi:\Sigma\mapsto X$ is $\mu_{\psi}$-almost
everywhere one-to-one and $\pi_{\ast}\mu_{\psi}= \mu_{\varphi}$.

Gibbs measures include all practically interesting invariant
measures on $X$ and $\Sigma$, e.g.~all smooth invariant measures
on $X$, Sinai-Ruelle-Bowen (SRB) measures, measures of maximal
entropy (i.e.~Margulis measures on $X$ and Parry measures on
$\Sigma$) etc.

We first study 
topological Markov
chains separately from Anosov diffeomorphisms. Let $\Sigma$
be a topological Markov chain with a transitive matrix $\bf A$.
Let $\mu$ be an arbitrary Gibbs measure defined by a H\"older
continuous potential. Naturally interesting subsets of $\Sigma$
are {\sl cylinders}, which include all balls in the metric $d_a$.

A cylinder $C\subset\Sigma$ is obtained by fixing symbols
on a finite interval $\Lambda=[n^-,n^+]\subset\ZZ$, i.e.~for some
$\omega_{\Lambda}\in\{1,\ldots,M\}^{\Lambda}$,
$\omega_{\Lambda}=\{\omega_{n^-},\ldots,\omega_{n^+}\}$, we set
\be
   C=C(\omega_{\Lambda});=\{\uom'\in\Sigma:\, \omega_i'=\omega_i\ \
         \ \ {\rm for}\ \ n^-\leq i\leq n^+\}
            \label{Comega}
\ee
Each cylinder is open and closed in $\Sigma$.
We call $n^-$ and $n^+$ the {\sl left} and {\sl right endpoints} of an
interval $\Lambda$, respectively, and $(n^-+n^+)/2$ the {\sl center} of
$\Lambda$.


Note that not every sequence of cylinders is a BC sequence.
For example, let $C_n=\sigma^nC$ for a fixed cylinder $C$.
It is obviously not a BC sequence.
Hence, we need some restrictions on cylinders
to ensure quasi-independence of $\sigma^{m-n}C_n$
and $C_m$ for large $|m-n|$.\medskip

\noindent{\bf Definition}. We say that two intervals
$[n_1^-,n_1^+]$ and $[n_2^-,n_2^+]$ are $D$-{\sl nested} for
$D\geq 0$ if either $[n_1^-,n_1^+]\subset [n_2^--D,n_2^++D]$ or
$[n_2^-,n_2^+]\subset [n_1^--D,n_1^++D]$.

\begin{theorem}
Let $\{C_n\}$ be a sequence of cylinders
defined on intervals $\Lambda_n\subset\ZZ$.
Let $D\geq 0$ be a constant. Assume that for all $m,n$ the
intervals $\Lambda_m,\Lambda_n$ are $D$-nested.
Then $\{C_n\}$ satisfies {\rm (SP)} and hence, if in addition $\sum \mu(C_n)
= \infty$,
it is an sBC sequence and  {\rm (\ref{log})} holds. \label{tm1}
\end{theorem}

\noindent
{\bf Examples}. \\
{\bf 1}. Let the left endpoints of $\Lambda_n$ lie in the interval
$[0,D]$, then $\Lambda_n$ are $D$-nested. We call such
intervals  $\Lambda_n$ 
$D$-{\sl aligned}. (Similarly one can talk about right
endpoints.) \\
{\bf 2}. Let the centers of $\Lambda_n$ lie in the interval
$[-D/2,D/2]$, then $\Lambda_n$ are $D$-nested.
We call such intervals $\Lambda_n$ 
$D$-{\sl centered}. Note that cylinders  defined on
$0$-centered intervals are precisely balls in $\Sigma$ with
respect to the metric $d_a$ defined above.
Therefore  the
``quantitative orbit density'' phenomenon (see the discussion after
Theorem~\ref{tmP}) holds for Gibbs measures on topological Markov
chains. Specifically, if  one fixes $\uom_0 \in \Sigma$ and considers
``a target shrinking
to $\uom_0$'', that is, a sequence  of balls (or centered at $\uom_0$,
then $\mu$-almost all orbits $\{\sigma^n\uom\}$ get into infinitely many such
balls whenever the sum of their measures diverges.

\medskip

The following two theorems show that the assumptions of
Theorem~\ref{tm1} cannot be easily relaxed.
We need to introduce some terminology generalizing the two examples
above.
Let $\{l_n\}$ be a sequence of positive numbers.
We say that  a sequence $\{\Lambda_n\}$ of intervals
is $\{l_n\}$-{\sl centered} (resp.~$\{l_n\}$-{\sl aligned}) if the
center (resp.~the left endpoint) of each $\Lambda_n$ belongs to
$[-l_n/2,l_n/2]$ (resp., $[0,l_n]$).

\begin{theorem}
Let $\{l_n\}$ be a sequence of natural numbers such that $l_n\to\infty$.
Then there is a sequence of cylinders $\{C_n\}$ with divergent sum of
measures
which is defined on $\{l_n\}$-centered (or,
alternatively, $\{l_n\}$-aligned) intervals $\Lambda_n\subset\ZZ$
 and does not satisfy {\rm (SP)}.
\label{tm2}
\end{theorem}

\begin{theorem}
Let $\varepsilon>0$.
There is a sequence of cylinders $\{C_n\}$ with $\sum \mu(C_n) =
\infty$ which is defined on $\{\varepsilon|\Lambda_n|\}$-centered (or,
alternatively, $\{\varepsilon|\Lambda_n|\}$-aligned) intervals
$\Lambda_n\subset\ZZ$
and is not a BC sequence. Moreover,
for a.e.~$\uom\in\Sigma$ there are only finitely many
$n$ such that $\sigma^n\uom\in C_n$.
\label{tm3}
\end{theorem}


Theorems~\ref{tm2} and \ref{tm3} show that it is not enough,
even for the BC property, that the cylinders are `relatively well'
centered or aligned.

\medskip

\noindent
{\bf Remarks}. \\
{\bf 1}. Suppose that  each of the
sets $C_n$ is a union of at most $k_n$ cylinders satisfying the
nested condition. It is clear that the conclusion of Theorem~\ref{tm1}
still holds  when the sequence $\{k_n\}$ is bounded.  On the other
hand, Theorem~\ref{tm2} shows that a sequence of unions $C_n$ of
$k_n$ $0$-centered cylinders may not satisfy (SP) if
$\{k_n\}$ is unbounded, while Theorem~\ref{tm3} shows that $\{C_n\}$
is not necessarily BC if $k_n$ is of order $n^a$ with some $a>0$. \\
{\bf 2}. Consider a one-sided topological Markov chain
$\sigma:\Sigma^+\mapsto\Sigma^+$ defined on the space $\Sigma^+$ of one-sided
sequences: 
$$
   \Sigma^+=\{\uom\in\{1,\ldots,M\}^{\ZZ_+}:\
    {\bf A}_{\omega_i\omega_{i+1}}=1\ \ \ \forall i\in\ZZ_+\}\,;
$$
here $\ZZ_+=\{0,1,2,\ldots\}$. Note that the shift $\sigma$
preserves $\Sigma^+$ but is not invertible, every sequence $\uom$
may have up to $M$ preimages. One-sided topological Markov chains
give symbolic representation for piecewise smooth expanding
interval maps satisfying the Markov condition.

Theorems~\ref{tm1}--\ref{tm3} apply to one-sided topologically
mixing Markov chains without change. Note, however, that all the
cylinders must be defined on intervals $\Lambda\subset\ZZ_+$. In
particular, our theorems hold for cylinders defined on intervals
that are $D$-aligned, $\{l_n\}$-aligned and
$\{\varepsilon|\Lambda_n|\}$-aligned, respectively\footnote{Note that in
this case the result of Theorem~\ref{tm1} can be derived from a
recent manuscript by Conze and Raugi \cite{CR}.}. Consider the metric
$d_a^+$ on $\Sigma^+$ given by $d_a^+(\uom,\uom')=a^n$ where
$n=\max\{n:\, \omega_i=\omega_i',\ \forall i<n\}$. In this metric,
balls are cylinders defined on $0$-aligned intervals. Therefore
the ``quantitative orbit density'' phenomenon, which follows from
Theorem~\ref{tmP} if $\Sigma^+=\{1,\ldots,M\}^{\ZZ_+}$ and $\mu$
is the product measure,  is extended to hold for an arbitrary
Gibbs measure on a one-sided topological Markov chain.

It is also worthwhile to mention that Theorem~\ref{tm3} gives
examples of non-BC sequences of cylinders in the setting of
one-sided shifts. In fact, the idea of the proof works for an
arbitrary measure-preserving system and produces examples of
non-BC sequences in the generality of Proposition~\ref{lmNBC}.

\medskip

Back to Anosov diffeomorphisms, the above theorems can be restated
by replacing cylinders with their projections on the manifold $X$
and the $T$-invariant measure $\pi_{\ast}\mu$ on $X$.
The projection $\pi(C)$ of a cylinder $C=C(\omega_{\Lambda})$
is a rectangle
\be
   \pi(C)=\cap_{i=n^-}^{n^+}T^{-i}R_{\omega_i}
    \label{piC}
\ee
in terms of of the formula (\ref{Comega}). These are very
special rectangles generated by the given Markov partition. It
would be of natural interest to extend our results to other
classes of rectangles, which we do next.

Recall that a rectangle $R$ is a subset of $X$ of a small diameter
such that for any points $x,y\in R$ the intersection $W^s_x\cap
W^u_y$ of the local stable manifold $W^s_x$ through $x$ and the
local unstable manifold $W^u_y$ through $y$ is a point that also
belongs in $R$. For $x\in R$ put $W^{u,s}_x(R)=W^{u,s}_x\cap R$.
For $x,y,\in R$ put $[x,y]=W^s_x\cap W^u_y$. Then for any $z\in R$
we have
$$
    R=[W^u_z(R),W^s_z(R)]=\{[x,y]:\, x\in W^u_z(R),y\in W^s_z(R)\}\,.
$$
So, $R$ has a direct product structure and
$W^u_z(R)$, $W^s_z(R)$ can be thought of as coordinate planes in
$R$. Note that $\partial R=\partial^uR\cup\partial^sR$,
where
$$
   \partial^uR=[W^u_z(R),\partial W^s_z(R)]
      \ \ \ \ \ \ {\rm and}\ \ \ \ \ \
   \partial^sR=[\partial W^u_z(R),W^s_z(R)]
$$
(these sets do not depend on $z\in R$).

We will consider small enough rectangles such that all local unstable
manifolds $W^u_x(R)$, $x\in R$ are almost parallel, and so are all
stable manifolds $W^s_x(R)$, $x\in R$. Hence, the diameters of our
rectangles are $\leq\varepsilon_1$ with some fixed small $\varepsilon_1>0$.
Our rectangles are not necessarily connected.

Our main assumption must be some sort of `roundness' of
rectangles, the necessity of which we explained above.
For any $\varepsilon>0$ put
$$
     W^u_z(R,\varepsilon):=\{x\in W^u_z(R):\,
     {\rm dist}(x,\partial W^u_z(R))<\varepsilon\}
$$
and
\be
           R^u_z(\varepsilon):=[W^u_z(R,\varepsilon),W^s_z]\,.
             \label{Ru}
\ee
This is a sort of  $\varepsilon$-neighborhood of the
stable boundary $\partial^sR$. Similarly, the $\varepsilon$-neighborhood
of the unstable boundary $\partial^uR$ is defined, call it
$R^s_z(\varepsilon)$.

Now fix another constant $\varepsilon_0 \in (0,\varepsilon_1)$
and some constants $C_0>0$, $\gamma>0$.
\medskip

\noindent{\bf Definition}. We say that a rectangle $R$ is
{\sl u-quasiround} if for some $z\in R$ \\
(i) the set $W^u_z(R)$ has (external) diameter $\leq\varepsilon_1$
and internal diameter $\geq\varepsilon_0$ (note that this set
will be perfectly round if $\varepsilon_0=\varepsilon_1$);
\\
(ii) For all $\varepsilon>0$
\be
   \mu( R^u_z(\varepsilon))\leq C_0|\ln\varepsilon|^{-1-\gamma}\mu(R)
        \label{lne}
\ee
Similarly, {\sl s-quasiround} rectangles are defined.
\medskip

Note that the definition of u- and s-quasiroundness
depends on the pre-fixed constants $\varepsilon_1,\varepsilon_0,
C_0,\gamma$.

The choice of $z$ in this definition is not important,
since the same properties will also holds for all $z\in R$,
with possibly slightly different values of
$\varepsilon_1,\varepsilon_0$ and $C_0$. The exact values of
$\varepsilon_1,\varepsilon_0,C_0,\gamma$ may affect some constants
in our estimates, but otherwise will be irrelevant.

Note that if the set $\partial W^u_z(R)$ is
smooth or piecewise smooth and the measure
on $W^u_z$ induced by $\mu$ is smooth, then
$\mu( R^u_z(\varepsilon))\leq {\rm const}\cdot\varepsilon\mu(R)$.
It is quite common in hyperbolic dynamics to assume that the
measure of $\varepsilon$-neighborhoods of boundaries
or singularities is bounded by const$\cdot\varepsilon^a$
for some $a>0$. Our bound (\ref{lne}) is milder than that.

Next, we need to consider arbitrary small rectangles that satisfy some
sort of roundness condition.
\medskip

\noindent{\bf Definition}. We call a rectangle $R$ {\sl
eventually quasiround} (EQR) if there are two integers $k^-\leq
k^+$ such that $T^{k^+}(R)$ is u-quasiround and $T^{k^-}(R)$ is
s-quasiround. \medskip

The integers $k^{\pm}$ may not be uniquely defined
for a rectangle $R$, but each of them is defined by $R$
up to a small additive
depending on the ratio $\varepsilon_1/\varepsilon_0$,
so the choice of $k^{\pm}$ for a given $R$ will not be important.

EQR rectangles in the Anosov setting play a role similar to that
of cylinders for TMC's, and the numbers $k^-$, $k^+$ correspond to
the endpoints of cylinders. Note, however, that EQR rectangles are
not generated by any Markov partitions. On the other hand, we
impose the regularity condition (\ref{lne}) on the boundary of EQR
rectangles, while no such condition was assumed for cylinders.

Note that if dim$\, X=2$, then stable and unstable manifolds are
one-dimensional, and, with appropriate choice of $\varepsilon_0$,
$\varepsilon_1$, every connected rectangle is EQR. Indeed, the
property
(i) follows from the uniform hyperbolicity of $T$ and the
compactness of $X$, while the property (ii) follows from our
Lemma~\ref{lmF2} in Section~\ref{secPAD} (note that the set
$R^u_z(\varepsilon)$ in this case consists of two connected
rectangles).\medskip

\noindent{\bf Definition}. We say that two EQR rectangles
$R_1,R_2$ with the corresponding integers $k_1^-,k_1^+$ and
$k_2^-,k_2^+$ characterizing their quasiroundness are $D$-{\sl
nested} for $D\geq 0$ if either $[k_1^-,k_1^+]\subset
[k_2^--D,k_2^++D]$ or $[k_2^-,k_2^+]\subset [k_1^--D,k_1^++D]$.
\medskip

\begin{theorem}
Let $T:X\mapsto X$ be an Anosov diffeomorphism with a Gibbs measure
$\mu$ defined by a H\"older continuous potential $\varphi$ on $X$,
and $D\geq 0$ a constant. Let $\{R_n\}$ be a sequence of EQR
rectangles. Assume that for all $m,n\geq
1$ the rectangles $R_m,R_n$ are $D$-nested. Then $\{R_n\}$
satisfies {\rm (SP)} and hence, if in addition $\sum \mu(R_n)
= \infty$,
it is an sBC sequence and  {\rm (\ref{log})} holds. \label{tm4}
\end{theorem}

\noindent{\bf Examples}.\\ {\bf 3}. If a sequence of EQR rectangles
$R_n$ satisfies the condition
\be
       |k^-_n+k^+_n|\leq D=\, {\rm const}
          \label{kn+-}
\ee
then it is an sBC sequence and verifies (\ref{log}).\\ {\bf 4}. In
particular, if $T$ is a linear 2-D toral automorphism and $\mu$
the Lebesgue measure, then any sequence of connected rectangles
with uniformly bounded ratio of stable and unstable sides (which
is sometimes called `aspect ratio') satisfies the condition
(\ref{kn+-}) and hence the conclusion of Theorem~\ref{tm4} holds. \\
{\bf 5}. Let 
$T:X\mapsto X$ be
the baker's transformation of the unit square $X=[0,1]\times
[0,1]$ and $\mu$ the Lebesgue measure. Note that $T$ is
discontinuous but still admits a finite Markov partition. Then any
sequence of balls with diverging measures is a BC sequence.
Indeed, in each ball $B\subset X$ one can find a `dyadic' square
$R\subset B$ such that $\mu(R)\geq 0.1\mu(B)$. Dyadic squares
correspond to 0-centered cylinders in the symbolic space, so one
can apply Theorem~\ref{tm1} and obtain the sBC property for the
dyadic squares, which implies (at least) the BC property for the
original balls.\medskip

Next, we generalize Example 4 to nonlinear Anosov diffeomorphisms.
Let $T:\, X\mapsto X$, dim$\, X=2$, be an Anosov diffeomorphism of a surface.
Recall that in this case every connected rectangle $R\subset X$ is EQR.
For a connected rectangle $R$ we denote
$$
     d^u(R)=\sup_{z\in R}|W^u_z(R)|
      \ \ \ \ \ \ {\rm and}\ \ \ \ \ \
     d^s(R)=\sup_{z\in R}|W^s_z(R)|\,,
$$
where $|W^u|$, $|W^s|$ stand for the Lebesgue measures
(lengths) of the corresponding curves $W^u, W^s$. Let $B\geq 1$.
We say that a rectangle $R$ has a $B$-{\sl bounded aspect ratio}
if
$$
      B^{-1}\leq d^u(R)/d^s(R)\leq B\,.
$$
Note that rectangles with $B$-bounded aspect ratio are, in the
geometric sense, close to squares (i.e., `round'). This geometric
version of roundness is somewhat more preferable and easier to
check than the dynamical roundness assumed by (\ref{kn+-}).

\begin{theorem}
Let $T:\, X\mapsto X$, dim$\, X=2$, be an Anosov diffeomorphism with
a Gibbs measure $\mu$ defined by a H\"older continuous potential
$\varphi$ on $X$, and $B\geq 1$ a constant. Let $\{R_n\}$ be a
sequence of connected rectangles with (uniformly) $B$-bounded aspect ratio.
Then $\{R_n\}$ satisfies {\rm (SP)} and hence, if in addition $\sum \mu(R_n)
= \infty$,
it is an sBC sequence and  {\rm (\ref{log})} holds. \label{tm5}
\end{theorem}

The extensions of Theorems~\ref{tm2} and \ref{tm3} to EQR
rectangles can also be obtained but are hardly worth pursuing,
because the examples of cylinders constructed in \ref{tm2}
and \ref{tm3} can be simply projected on $X$ and produce the
corresponding examples of rectangles.

\section{Proofs for topological Markov chains}
\label{secPTMC}
\setcounter{equation}{0}

The following facts about Gibbs measures are standard:

\noindent
{\bf Fact 1} \  For any cylinder $C$ defined on an interval $\Lambda$
$$
     c_1\theta_1^{|\Lambda|}\leq\mu(C)\leq c_2\theta_2^{|\Lambda|}\,,
$$
where $c_1,c_2>0$ and $\theta_1,\theta_2\in (0,1)$ only depend
on the Gibbs measure $\mu$.

\noindent
{\bf Fact 2} \ Let $C_1\subset C$ be cylinders
defined on intervals $\Lambda_1,\Lambda$ (note that
in this case $\Lambda_1\supset\Lambda$), then
$$
     c_1\theta_1^{|\Lambda_1|-|\Lambda|}
      \leq\mu(C_1)/\mu(C)\leq c_2\theta_2^{|\Lambda_1|-|\Lambda|}\,.
$$

\noindent
{\bf Fact 3} \ Let $C_1,C_2$ be cylinders defined on disjoint
intervals $[n_1^-,n_1^+]$ and $[n_2^-,n_2^+]$ in $\ZZ$.
Assume, without loss of generality that $n_1^+<n_2^-$. Then
$$
    |\mu(C_1\cap C_2)-\mu(C_1)\mu(C_2)|\leq
      c_3\theta_3^{n_2^--n_1^+}\mu(C_1)\mu(C_2)\,,
$$
where $c_3>0$ and $\theta_3\in (0,1)$ only depend on the Gibbs
measure $\mu$.\medskip

Facts 1 and 2 can be proved with the help of a normalized
potential for the Gibbs measure $\mu$, see \cite{ChS}. Fact 3 is
proved by R.~Bowen in \cite{Bo75}.

Let us introduce the following notation. If
$\Lambda_1=[n_1^-,n_1^+]$ and $\Lambda_2=[n_2^-,n_2^+]$ are two
intervals (not necessarily disjoint), define an ``asymmetric
distance'' $\delta(\Lambda_1,\Lambda_2)$ by
$$
 \delta(\Lambda_1,\Lambda_2)=\min\{D:\, \Lambda_2{\rm\ is\ in\
 the\ }D{\rm-neighborhood\ of\ }\Lambda_1\}
$$
Equivalently,
$\delta(\Lambda_1,\Lambda_2)=\max\{n_2^+-n_1^+,n_1^--n_2^-,0\}$.
Clearly, $\delta(\Lambda_1,\Lambda_2)=0$ if and only if
$\Lambda_2\subset\Lambda_1$. It is also clear that
$\Lambda_1,\Lambda_2$ are $D$-nested if and only if one of
the distances $\delta(\Lambda_1,\Lambda_2)$ and $\delta(\Lambda_2,\Lambda_1)$
does not exceed $D$. 

\begin{lemma}
If $C_1$, $C_2$ are cylinders defined on intervals $\Lambda_1$ and
$\Lambda_2$, respectively, then
$$
    |\mu(C_1\cap C_2)-\mu(C_1)\mu(C_2)|\leq
      c_4\theta_4^{\delta(\Lambda_1,\Lambda_2)}\mu(C_1)\,,
$$ where $c_4>0$ and $\theta_4\in (0,1)$ only depend on the Gibbs
measure $\mu$.
\label{lmCC}
\end{lemma}

{\em Proof}. This follows from Facts 1 and 3 if $\Lambda_1$ and
$\Lambda_2$ are disjoint, and from Facts 1 and 2 if they are not.
$\Box$\medskip

{\em Proof of Theorem~\ref{tm1}}. We estimate the quantity
$R_{mn}=\mu(C_m\cap\sigma^{m-n}C_n)-\mu(C_n)\mu(C_m)$. Without
loss of generality, assume that the interval $\Lambda_m$ is
``nested'' in $\Lambda_n$, i.e.~$\Lambda_m$ lies in the
$D$-neighborhood of $\Lambda_n$.
Note that we do not assume any
relation between $m$ and $n$, or  between $\mu(C_m)$ and
$\mu(C_n)$. Our assumption easily implies that $\Lambda_n-(m-n)$
is not in the $(|m-n|-D)$-neighborhood of $\Lambda_m$. Applying
Lemma~\ref{lmCC} to the cylinders $C_m$ and $\sigma^{m-n}C_n$, one
gets
$$
   |R_{mn}|\leq c_4\theta_4^{\delta(\Lambda_m,\Lambda_n-(m-n))}\mu(C_m)
   \leq c_4\theta_4^{|m-n|-D}\mu(C_m)
$$
Summing up over all $n$ satisfying our nesting condition (that
$\Lambda_m$ is ``nested'' in $\Lambda_n$) gives a quantity bounded
by const$\cdot\mu(C_m)$. Now summing up over $m=M,\ldots,N$ proves
(SP). $\Box$\medskip

In the following proofs of Theorems~\ref{tm2} and \ref{tm3} we use
a special construction. Let $T$ be a measure preserving transformation
(invertible or not) of a probability space $(X,\mu)$, and let $\{\tilde{A}_k\}$ be a sequence of
measurable subsets of $X$  and $\{l_k\}$ a sequence of natural
numbers. Put $s_0=0$
and $s_k=l_1+\cdots +l_k$ for $k\geq 1$. Consider a new sequence
of sets $\{{A}_n\}$ defined as follows:
$$
   T^{1-l_1}\tilde{A}_1,T^{2-l_1}\tilde{A}_1,\ldots,
              T^{-1}\tilde{A}_1,\tilde{A}_1,
   T^{1-l_2}\tilde{A}_2,\ldots,T^{-1}\tilde{A}_2,\tilde{A}_2,
   T^{1-l_3}\tilde{A}_3,\ldots,T^{-1}\tilde{A}_3,\tilde{A}_3,\ldots
$$
Note that the $n$th set in this sequence is
\be
       A_n=T^{n-s_k}\tilde{A}_k\,,
        \label{CnCk}
\ee
where $k$ is defined by $s_{k-1}<n\leq s_k$. We denote this
$k$ by by $k=k_n$. We will say that the new sequence, $\{A_n\}$,
is {\sl derived} from $\{\tilde{A}_k\}$ and $\{l_k\}$. \medskip

{\em Proof of Theorem~\ref{tm2}}. Without loss of generality,
assume that $l_n$ is monotonic, $1\leq l_1\leq l_2\leq\cdots$.
Let $\{\tilde{C}\}$ be a cylinder defined on some
interval $[0,l]$
(alternatively, we can assume that its center is at zero).
Now consider the sequence of cylinders $\{C_n\}$ derived from
the constant sequence $\tilde{C}_k = \tilde{C}$ and $\{l_k\}$. Then
$C_n$ is defined on
an interval $\Lambda_n$ whose left endpoint lies
in the interval $[0,l_k]$ where $k=k_n$ is defined above. Since
$\{l_k\}$ is monotonic, the left endpoint of $\Lambda_n$ lies in $[0,l_n]$,
so the nesting condition of Theorem~\ref{tm2} is satisfied. It is
now easy to see that for $N=l_1+\cdots+l_k$ we have $
   E_N =  \sum_{n = 1}^N \mu(C_n)= (l_1 + \cdots + l_k)
\mu(\tilde{C})$, while
$$
     \sum_{m,n=1}^NR_{mn}\geq
     \frac 12\,\Big (l_1^2+\cdots +l_k^2\Big )\,\mu(\tilde{C})\,.
$$
It is clear that the right hand side of this inequality grows
faster than $CE_N$ for any $C>0$, which violates (SP). $\Box$\medskip

We write $a_n\approx b_n$ for two sequences of numbers $\{a_n\}$ and
$\{b_n\}$ if there are constants $0<c_1<c_2<\infty$ such that
$c_1<a_n/b_n<c_2$ for all $n$ (the constants $c_1,c_2$ may depend
on the topological Markov chain $(\Sigma_{\bf A},\sigma)$ and the
Gibbs measure $\mu$).\medskip

{\em Proof of Theorem~\ref{tm3}}. Let $\{\tilde{C}_k\}$ be a
sequence of cylinders defined on
intervals $\tilde{\Lambda}_k$
with left endpoints at zero such that
$\mu(\tilde{C}_k)\approx 1/(k\ln^2k)$. (Again, we could assume
that the centers of $\{\tilde{\Lambda}_k\}$ are at zero.)
It follows from Fact 1 that
$|\tilde{\Lambda}_k|\approx\log k$. For each $k\geq 1$, let
$l_k=[\varepsilon |\tilde{\Lambda}_k|]$. Consider the sequence of
cylinders $\{C_n\}$ derived from $\{\tilde{C}_k\}$ and $\{l_k\}$.
Then $C_n$ is defined on an interval $\Lambda_n$ whose
left endpoint lies in the interval $[0,\varepsilon |\Lambda_n|]$.
Since $l_k\approx\log k$, we have
$\sum\mu(C_n)=\sum_k l_k\mu(\tilde{C}_k)=\infty$. On the other
hand, $\sum_k\mu(\tilde{C}_k)<\infty$. Hence, by Lemma~\ref{lmBCT}
(i), for a.e.~$\uom\in \Sigma$ there are at most finitely many $k\geq
0$ such that
$\sigma^{s_k}\uom\in \tilde{C}_k$. Now, by (\ref{CnCk}),
$C_n=\sigma^{n-s_k}\tilde{C}_k$, hence there are
 at most finitely many $n$ such that $\sigma^n\uom\in C_n$. $\Box$\medskip

Lastly, we give proofs of  two propositions from the
introduction.

{\em Proof of Proposition~\ref{lmWM}}. Part (i) easily follows from the
ergodic theorem. For part (ii), let $T$ be weakly mixing and
$\{A_n\}$ contain
finitely many distinct subsets of $X$ of positive measure, call them
$F_1,\ldots,F_k$.
Since $c_1<E_N/N<c_2$ for some constants $0<c_1<c_2<\infty$, to show
that
\be
\mu(S_N/E_N-1)^2\to 0\ \ \ \ \ {\rm as }\ \ N\to\infty\label{L2}
\ee
it is enough to prove that
\be
    \mu(S_N-E_N)^2=\sum_{m,n=1}^NR_{mn}=o(N^2)
    \label{ON2}
\ee
The weak mixing of $T$ implies that for any $F_i,F_j$
$$\sum_{n=1}^N|\mu(T^{-n}F_i\cap F_j)-\mu(F_i)\mu(F_j)|=o(N)\,,$$
and since we only have finitely many pairs $(F_i,F_j)$, the term
$o(N)$ here is uniform in $i,j$. This completes the proof of
(\ref{ON2}).  On the other hand, if (\ref{L2}) holds,  one can  choose
a subsequence $\{N_k\}$ such that
$S_{N_k}/E_{N_k}\to 1$ almost surely. Thus $S_{N_k}\to\infty$ on
a set of full measure, which clearly implies
that $\{A_n\}$ is a BC sequence.

Assume now that $T$ is not weakly mixing. If it is not ergodic,
the constant sequence $A_n = A$, where $A$ is a nontrivial
invariant set, is clearly not a BC sequence. Otherwise $T$
has a factor isomorphic to a rotation of a circle (because $T$ has
a non-constant eigenfunction with eigenvalue $\exp(2\pi\theta i)$
with some $0 < \theta <1$, see e.g.~\cite{Pet}, p. 65--68). If
$\theta$ is rational, then $T^k$ is not ergodic for some $k$ and
the claim follows as above. If $\theta$ is irrational, then the
factor measure is Lebesgue. To finish the proof of (ii) it is then
enough to consider an
irrational rotation of a circle and find a sequence of (nonempty)
arcs $\{A_n\}$ that only contains finitely many distinct arcs but
is not BC. This is a simple exercise. Part (iii) follows from the
definitions in a straightforward way and is also left as an exercise to
the reader. $\Box$\medskip

{\em Proof of Proposition~\ref{lmNBC}}. If for some $\varepsilon$ there are
no measurable
subsets $A$ of $X$ with $0 < \mu(A) < \varepsilon$, then (assuming
$\mu$ is nontrivial)  $T$ is
not weakly mixing, so the claim follows from the previous proposition.
Otherwise there exists  a sequence $\{\tilde{A}_k\}$ of sets of
positive measure such that $\sum_{k = 1}^\infty \mu(\tilde{A}_k)<
\infty$. Define a sequence  $\{l_k\}$  of natural numbers by $l_k
= [1/\mu(\tilde{A}_k)] + 1$, and let $\{{A}_n\}$ be a sequence
{derived} from $\{\tilde{A}_k\}$ and $\{l_k\}$. Then clearly
$\sum_{n = 1}^\infty \mu({A}_n)= \infty$.  On the other hand, we
can argue as in the proof of Theorem~\ref{tm3} to show that  for
a.e.~$x\in X$  there are
 at most finitely many $n$ such that $T^n x\in A_n$. $\Box$\medskip

\section{Proofs for Anosov diffeomorphisms}
\label{secPAD}
\setcounter{equation}{0}

In this section we use an approach based on the shadowing property and
specification. Ruelle recently demonstrated the power and elegance
of this approach in \cite{Ru99}, and we follow his lines.

We recall certain standard facts about transitive Anosov
diffeomorphisms. We will denote by $\Lambda$ finite or infinite
intervals of $\ZZ$. For a finite interval $\Lambda=[n^-,n^+]$, we
denote by $|\Lambda|=n^+-n^-+1$ the cardinality of $\Lambda$. For
two disjoint intervals $\Lambda_1,\Lambda_2$ we denote by
$$
 {\rm dist}(\Lambda_1,\Lambda_2)=\min\{|i-j|:\, i\in\Lambda_1,j\in\Lambda_2\}
$$
the length of the gap between them.
\medskip

\noindent
{\bf Expansiveness}. Any Anosov diffeomorphism $T:X\mapsto X$
is {\sl expansive}, i.e.~there is a $\delta>0$ (called {\sl expansivity
constant}) such that
$$
       \forall k\in\ZZ\ \ d(T^kx,T^ky)<\delta\ \ \ \
         \Leftrightarrow\ \ \ \ x=y\,.
$$
In fact, due to the hyperbolicity of $T$, for some $C>0$ and  $0 <
\theta < 1$ one has
\be
       \forall |k|\leq n\ \ d(T^kx,T^ky)<\delta\ \ \ \
         \Rightarrow\ \ \ \ d(x,y)<C\theta^{-n}\,.
             \label{exp1}
\ee

Let $\Lambda\subset\ZZ$ be an interval of $\ZZ$, finite or not.
Let ${\bf x}=(x_k)_{k\in\Lambda}\in X^{\Lambda}$. Given $\alpha>0$,
we say that $\bf x$ is an $\alpha$-{\sl pseudo-orbit} if
$$
     d(T^kx,x_{k+1})<\alpha\ \ \ \ {\rm whenever}\ \ \ \ k,k+1\in\Lambda\,.
$$
We say that the orbit of $x\in X$ $\beta$-{\sl shadows} $\bf x$ if
$$
      d(T^kx,x_k)<\beta\ \ \ \ \ \ \forall k\in\Lambda\,.
$$
\medskip

\noindent
{\bf Shadowing lemma}. For any $\beta>0$ there is an $\alpha>0$
such that every $\alpha$-pseudoorbit is $\beta$-shadowed by a true
orbit of some $x\in X$.

Note that if $\Lambda=\ZZ$ and $\beta<\delta/2$,
then the true orbit shadowing $\bf x$ is unique
by the expansivity. We fix a
$\beta<\delta/2$ and this fixes the corresponding
$\alpha>0$.

Note that if the pseudoorbit is periodic, then it is shadowed
by a true periodic orbit with the same period.

Given $\alpha>0$, there is an integer $K>0$ such that for every
$x,y\in X$ and $n\geq K$ there is a $z\in X$ such that
$$
     d(z,x)<\alpha\ \ \ \ \ {\rm and}\ \ \ \ \ d(T^nz,y)<\alpha\,,
$$
which follows from the topological transitivity of $T$.
(Note that our choice of $\alpha$ made above also fixes $K$.)

Using this remark, we can interpolate (concatenate)
several $\alpha$-pseudoorbits defined on intervals
of $\ZZ$ separated by gaps of lengths $\geq K$
in the following way.

\medskip
\noindent{\bf Specification}. Let $\alpha$-pseudoorbits
${\bf x}_j$ be defined on disjoint intervals of $\ZZ$ separated by gaps
of length $\geq K$. Then the ${\bf x}_j$ are all $\beta$-shadowed
by one true orbit of some $x\in X$.

One can also find a periodic orbit that $\beta$-shadows all
${\bf x}_j$, with period $P:=i_{\max}-i_{\min}+K$, where $i_{\max}$
and $i_{\min}$ are the maximum and the minimum points of the
union of the intervals of $\ZZ$ on which the pseudoorbits
${\bf x}_j$ are defined.

Due to the expansivity, the number of periodic orbits of period $P$
in the above construction is less than some $L$ independent of
the lengths of the intervals of $\ZZ$ where the pseudoorbits are
defined. The value of $L$ only depends on the number of these
intervals and the lengths of gaps between them. In our further
arguments, we will interpolate no more than four pseudoorbits
at a time, and the gaps between them will never exceed $2K$,
so we just fix the corresponding constant $L$.

Now, let $g:X\mapsto \IR$ be a H\"older continuous function. The
bound (\ref{exp1}) implies the following.

\medskip
\noindent{\bf Approximation of sums along orbits}.
There is a constant $B=B(g)$ such that
$$
   \forall k\in [p,q]\ \ \ d(T^kx,T^ky)<\delta\ \ \ \ \
      \Rightarrow\ \ \ \ \ \left |\sum_{k=p}^qg(T^kx)
      -\sum_{k=p}^qg(T^ky)\right |\leq B\,.
$$
Furthermore, let the specification property be used to shadow
two finite orbits $\{T^kx'\}$, $k\in\Lambda'$, and
$\{T^kx''\}$, $k\in\Lambda''$, with
$$
   K\leq\, {\rm dist}(\Lambda',\Lambda'')\leq 2K\,,
$$
by a periodic orbit of $z$ of period
$$
    P=|\Lambda'|+|\Lambda''|+\,{\rm dist}(\Lambda',\Lambda'')+K\,,
$$
then
\be
      \left |\sum_{k\in\Lambda'}g(T^kx')+\sum_{k\in\Lambda''}g(T^kx'')
      -\sum_{k=1}^{P}g(T^kz)\right |\leq B':=2B+3K||g||_{\infty}\,.
         \label{3B}
\ee
Note that $B'$ is a constant, just like $B$, independent of the
lengths of the intervals $\Lambda',\Lambda''$.

For $n\geq 1$, let
$$
     {\rm Fix}(T^n,X)=\{x\in X:\, T^nx=x\}
$$
be the set of periodic points of period $n$ in $X$.

\medskip
\noindent{\bf Periodic orbit approximation of Gibbs measures}.
Let $\mu$ be a Gibbs measure corresponding to a H\"older
continuous potential $\varphi:X\mapsto \IR$. For each $n\geq 1$, let
$\mu_n$ be an atomic probability measure concentrated on
Fix$(T^n,X)$ that assigns weight
\be
      \mu_n(x)=Z_n^{-1}\exp[\varphi(x)+\varphi(Tx)+\cdots+\varphi(T^{n-1}x)]
           \label{mun}
\ee
to each point $x\in\,$Fix$(T^n,X)$ (here $Z_n$ is a normalizing factor).
Then $\mu_n$ weakly converges to $\mu$ as $n\to\infty$.

\medskip
\noindent{\bf Variational principle}.
Let $\varphi:X\mapsto\IR$ be a continuous function and $P_\varphi$ its
topological pressure. Then
\be
      \sup_{\nu}[h_{\nu}(T)+\nu(\varphi)]=P_\varphi\,,
           \label{varpr}
\ee
where the supremum is taken over all $T$-invariant probability
measures $\nu$ on $X$, and $h_{\nu}(T)$ is the Kolmogorov-Sinai
entropy of $\nu$. Any measure $\nu$ that turns (\ref{varpr}) into
an equality is called an equilibrium state for $\varphi$.
Equilibrium states exist for every continuous function $\varphi$.
If $\varphi$ is
H\"older continuous on $X$, the equilibrium state is unique and
coincides with the Gibbs measure for the potential $\varphi$.

We now prove a few technical lemmas. Let $\mu$ be a Gibbs measure
on $X$ corresponding to a H\"older continuous potential $\varphi$.

We generalize our notation of Section~\ref{secPTMC}
by writing for any two variable quantities  $A$ and $B$
$$
     A\approx B\ \ \ \ \ \Leftrightarrow \ \ \ \ \ 0<c_1<A/B<c_2<\infty
$$
for some constants $c_1,c_2$ that only depend on $T:X\mapsto X$ and
the Gibbs measure $\mu$.

\begin{lemma}
The normalizing factor (the analogue of partition function) $Z_n$
in {\rm (\ref{mun})} satisfies
$$
       Z_n\approx e^{P_\varphi n}\,.
$$
\label{lmZn}
\end{lemma}

Note that it is standard to compute the topological pressure as
$$
      P_\varphi=\lim_{n\to\infty}\frac 1n \ln Z_n\,.
$$
The estimate in our lemma is sharper than this standard formula.

We need an elementary sublemma that is a modification
of a standard one, see Lemma~1.18 in \cite{Bo75}.

\begin{sublemma}
Let $\{a_n\}_{n=1}^{\infty}$ be a sequence of real numbers
such that $|a_{m+n}-a_m-a_n|\leq R$ for all $m,n\geq 1$
and some constant $R>0$.
Then $P:=\lim_{n\to\infty}a_n/n$ exists.
Furthermore, $|a_n-Pn|\leq 2R$ for all $n$.
\end{sublemma}

{\em Proof}.
Fix an $m\geq 1$. For $n\geq 1$, write $n=km+l$
with $0\leq l\leq m-1$. Then it follows by induction
on $k$ that $|a_n-ka_m-a_l|\leq kR$. Hence,
$$
   \left |\frac{a_n}{n}-\frac{ka_m}{km+l}-\frac{a_l}{km+l}\right |
    \leq \frac{kR}{km+l}\,.
$$
Letting $n\to\infty$ gives
$$
   \frac{a_m}{m}-\frac Rm \leq
   \liminf_n\frac{a_n}{n}\leq
   \limsup_n\frac{a_n}{n}\leq
   \frac{a_m}{m}+\frac Rm\,.
$$
Hence, $P:=\lim a_n/n$ exists. Next, assume
that $a_m>Pm+2R$ for some $m$. Then
$a_{2^nm}>2^nmP+(2^n+1)R$
which follows by induction on $n$.
Hence $\lim\sup a_n/n\geq P+R/m$,
a contradiction. A similar contradiction results
from the assumption $a_m<Pm-2R$. $\Box$\medskip

{\em Proof of Lemma~\ref{lmZn}}.
It is enough to show that
$$
R:=\sup_{m,n}|\ln Z_{m+n}-\ln Z_{m}-\ln Z_{n}|<\infty
$$
and apply the previous sublemma to the sequence $a_n=\ln Z_{n}$.
So, we need to show that
$$
  Z_{m+n}\approx Z_{m} Z_{n}\,.
$$
For fixed $n,m$, put $\Lambda'=[0,m-K]$ and
$\Lambda''=[m,m+n-K]$. For any $x\in\,$Fix$(T^{m+n},X)$ consider
${\bf x}'=\{x,\ldots, T^{m-K}x\}$ and ${\bf x}''=\{T^{m}x,\ldots,
T^{m+n-K}x\}$, these are two pseudoorbits defined on the intervals
$\Lambda'$ and $\Lambda''$ separated by a gap of length $K$. Each
of them can be shadowed by a true periodic orbit, of periods $m$
and $n$, respectively, and there are at most $L$ of those periodic
orbits for each of ${\bf x}'$ and ${\bf x}''$. On the other hand,
for every pair of periodic orbits $y'\in\,$Fix$(T^m,X)$ and
$y''\in\,$Fix$(T^n,X)$ consider two pseudoorbits ${\bf
y}'=\{y',\ldots,T^{m-K}y'\}$ defined on $\Lambda'$ and ${\bf
y}''=\{y'',\ldots,T^{n-K}y'\}$ defined on the interval $\Lambda''$
by associating $T^iy''$ to $m+i\in\Lambda''$. Then there is a true
periodic orbit of period $m+n$ shadowing both ${\bf y}'$ and ${\bf
y}''$, and the number of those periodic orbits does not exceed
$L$. Now the result follows from (\ref{3B}). $\Box$\medskip

Note that the potential $\varphi-P_\varphi$
corresponds to the same measure $\mu$
and has zero topological pressure. Hence we may just assume that
$P_\varphi=0$ in what follows. Then $Z_n\approx 1$.

We assume, as we may, that $\varepsilon_1$ in the definition
of EQR rectangles does not exceed the expansivity constant
$\delta$.

\begin{lemma}
Let $R$ be an EQR rectangle with integers $k^+$ and $k^-$ characterizing
the quasiroundedness of $R$. Let $x\in R$. Then
$$
     \mu(R)\approx \exp[\varphi(T^{k^-}x)+\cdots+\varphi(T^{k^+}x)]\,.
$$
\label{lmR1}
\end{lemma}

{\em Proof}. Consider a pseudoorbit ${\bf
x}=\{T^{-k^-}x,\ldots,T^{k^+}x\}$ defined on the interval \linebreak 
$\Lambda=[-k^-,k^+]$. Let $n\gg|\Lambda|$. Note that $y\in R$ if
and only if (i) $W^s(T^{k^+}y)$ intersects $W^u(T^{k^+}x)$, and
(ii) $W^u(T^{k^-}y)$ intersects $W^s(T^{k^-}x)$. This definitely
happens if the orbit ${\bf y}=\{T^{k^-}y,\ldots,T^{k^+}y\}$
$\varepsilon_0$-shadows $\bf x$. On the other hand, if $y\in R$,
then $\bf y$ $\varepsilon_1$-shadows $\bf x$. Hence, we can apply
our previous estimates with $\alpha=\varepsilon_0$ and
$\alpha=\varepsilon_1$, the value of $\alpha$ only affects the
values of all constants, which are not essential. So, we may simply
assume that $y\in R$ if and only if $\bf y$ $\alpha$-shadows $\bf
x$.

Now, for any $y\in\,$Fix$(T^n,X)$ consider the pseudoorbit
${\bf y}=\{T^{k^++K}y,\ldots,T^{n+k^--K}y\}$ defined on
the interval $\Lambda'=[k^++K,n+k^--K]$. Then ${\bf y}$ is
shadowed by a true periodic orbit of period $p:=n-(k^+-k^-)-K+1$,
and the number of those periodic orbits is less than $L$.
On the other hand, for any $z\in\,$Fix$(T^p,X)$ consider
a pseudoorbit ${\bf z}=\{z,\ldots,T^{p-1}z\}$ defined on the interval
$\Lambda'$ by associating $T^iz$ to $k^++K+i\in\Lambda'$.
Then there is a true periodic orbit of period $n$ shadowing
both ${\bf x}$ and ${\bf z}$, and the number of those periodic
orbits does not exceed $L$. Now the result
follows from (\ref{3B}) and the facts $Z_n\approx 1$ and
$Z_p\approx 1$. $\Box$\medskip

\begin{lemma}
Let $R_1,R_2$ be two EQR rectangles with integers
$k^{\pm}_1$ and $k^{\pm}_2$ characterizing the quasiroundedness
of $R_1,R_2$. Denote $\Lambda_i=[k_i^-,k^+_i]$ for $i=1,2$.
Let $x\in R_1\cap R_2$. Then
$$
    \mu(R_1\cap R_2)\leq c\cdot\exp
      \left [\sum_{i\in\Lambda_1\cup\Lambda_2}\varphi(T^ix)\right ]\,,
$$
with some $c>0$ that only depends on the Gibbs measure $\mu$.
\label{lmR2}
\end{lemma}

{\em Proof}. The proof of the previous lemma applies with the
following simple adjustments. Note that if $y\in R_1\cap R_2$,
then the orbit of $y$ $\varepsilon_1$-shadows that of $x$ on
$\Lambda_1\cup\Lambda_2$. So, to get an upper bound on
$\mu(R_1\cap R_2)$, we can take into account all $n$-periodic
orbits that $\alpha$-shadow the orbit of $x$ on
$\Lambda_1\cup\Lambda_2$ with $\alpha=\varepsilon_1$. Now, if
$\Lambda_1$ and $\Lambda_2$ overlap, the argument is exactly like
in the proof of the previous lemma. Let $\Lambda_1$ and
$\Lambda_2$ be disjoint with dist$(\Lambda_1,\Lambda_2)=J$. If
$J\leq 2K$, we can simply disregard such a small gap and apply the
previous argument. If $J>2K$, we replace the part of the orbit of
$y\in\,$Fix$(T^n,X)$ of length $J$ between $\Lambda_1$ and
$\Lambda_2$ by periodic orbits of period $J-2K$. To conclude the
argument, we now need an obvious extension of (\ref{3B}) from two
to four pseudoorbits with gaps of length $K$ in between. This
extension is straightforward. $\Box$\medskip

\begin{lemma}
There is a constant $\Delta>0$ such that for all $n\geq 1$
and $x\in\,{\rm Fix}(T^n,X)$ we have
$$
   \varphi(x)+\varphi(Tx)+\cdots+\varphi(T^{n-1}x)\leq -\Delta n\,.
$$
\end{lemma}

{\em Proof}. Let $\delta_x$ be the delta measure concentrated
at $x$. The measure
$$
    \delta_{x,n}=\frac 1n (\delta_x+\cdots+\delta_{T^{n-1}x})
$$
is $T$-invariant, so by the variational principle we have
$$
    \delta_{x,n}(\varphi)=\frac 1n \Big(
\varphi(x)+\cdots+\varphi(T^{n-1}x)\Big )\leq 0\,.
$$
We now need to prove that
$$
     \sup_{n\geq 1}\sup_{x\in\,{\rm Fix}(T^n,X)}\delta_{x,n}(\varphi)<0\,.
$$
If this is not true,
then there is a sequence of periodic points $x_k\in\,$Fix$(T^{n_k},X)$
such that $\delta_{x_k,n_k}(\varphi)\to 0$. We take any limit point of the
sequence of measures $\delta_{x_k,n_k}$ in the weak topology,
it will be a $T$-invariant
measure, call it $\nu$. We have $\nu(\varphi)=0$, so by the uniqueness
part of the variational principle $\nu=\mu$, so $\mu(\varphi)=0$ and
hence $h_{\mu}(T)=0$.
But it is known that $h_{\mu}(T)>0$ for any Gibbs measure,
a contradiction. $\Box$\medskip

Combining this lemma with the specification property and (\ref{3B}) gives

\begin{corollary}
There is a constant $B_0>0$ such that for all
$n\geq 1$ and $x\in X$
$$
          -\Delta_0 n\leq
      \varphi(x)+\varphi(Tx)+\cdots+\varphi(T^{n-1}x)\leq B_0-\Delta n
$$
with $\Delta_0=||\varphi||_{\infty}$.
\label{crh}
\end{corollary}

We can now prove analogues of Facts~1 and 2 of
Section~\ref{secPTMC} for Anosov diffeomorphisms. Our constants,
such as $c_i,\theta_i$, will only depend on the Gibbs measure
$\mu$ and the values of $\varepsilon_0,\varepsilon_1$ in the
definition of EQR rectangles. We use notation of Lemmas~\ref{lmR1}
and \ref{lmR2}.

\begin{lemma}
Let $R$ be an EQR rectangle and $k:=k^+- k^-$. Then
$$
    c_5\theta_5^{k}\leq\mu(R)\leq c_6\theta_6^{k}
$$
with some $c_5,c_6>0$ and $\theta_5,\theta_6\in (0,1)$.
\label{lmF1}
\end{lemma}

\begin{lemma}
Let $R_1,R_2$ be EQR rectangles with the corresponding intervals
$\Lambda_i=[k_i^-,k^+_i]$, $i=1,2$. Put $k:=|\Lambda_2 \setminus
\Lambda_1|$. Then
$$
    \mu(R_1\cap R_2)\leq c_7\theta_7^{k}\mu(R_1)
$$
with some $c_7>0$ and $\theta_7\in (0,1)$.
\label{lmF2}
\end{lemma}

{\em Proof}. Lemmas~\ref{lmF1} and \ref{lmF2}
follow from Lemmas~\ref{lmR1} and \ref{lmR2}
and Corollary~\ref{crh}. $\Box$\medskip

Note that so far we only used the property (i) of the
quasiround rectangles, we did not use (\ref{lne}).

\begin{lemma}
Let $R_1,R_2$ be EQR rectangles such that the intervals
$\Lambda_1=[k^-_1,k^+_1]$ and $\Lambda_2=[k^-_2,k^+_2]$ are
disjoint. Put $k:=\,{\rm dist}(\Lambda_1,\Lambda_2)$. Then
$$
    |\mu(R_1\cap R_2)-\mu(R_1)\mu(R_2)|\leq
      c_8\frac{\mu(R_1)+\mu(R_2)}{|ak+b|^{1+\gamma}}
$$
with some constants $c_8>0$, $a>0$ and $b$.

\label{lmF3}
\end{lemma}

{\em Proof}. Our proof uses Markov partitions and symbolic
dynamics. Let ${\cal R}$ be a Markov partition and $\Sigma$ the
corresponding symbolic space, a topological Markov chain. We now
partition the rectangles $R_1$ and $R_2$ into subrectangles
generated by the Markov partition $\cal R$ as follows. Let
$C\subset\Sigma$ be a cylinder defined on an interval
$\Lambda\subset\ZZ$. We say that its projection $\pi(C)$ is {\sl
properly inside} $R_i$, $i=1,2$, if

 (i) $\pi(C)\subset R_i$,
and

 (ii) for any larger cylinder $C'\supset C$ its projection
$\pi(C')$ is not a subset of $R_i$.\\ Denote by ${\cal C}_i$ the
collection (in general, countable) of cylinders that are properly
inside $R_i$. Since $R_i$ is a rectangle, one can easily check
that all the cylinders in ${\cal C}_i$ are disjoint. Next, it
follows from the assumption (\ref{lne}) that $\mu(\partial
R_i)=0$, hence $$
    \mu\big (R_i\setminus\bigcup_{C\in{\cal C}_i}\pi(C)\big )=0\,,
$$
i.e.~the rectangles $\pi(C)$, $C\in{\cal C}_i$,
make a (mod 0) partition of $R_i$.

Now consider the collection ${\cal C}_1$ and an arbitrary cylinder
$C\in{\cal C}_1$ defined on an interval $\Lambda=[k^-,k^+]$.
Observe that if $t:=k^+-k_1^+>0$, then, using the notation of
(\ref{Ru}), we have $\pi(C)\subset R^u_{1,z}(\varepsilon)$ with
$\varepsilon=c\theta^t$ for any $z\in R_1$. Here $c>0$ and
$\theta\in (0,1)$ are constants determined by the hyperbolicity
properties of $T$ and the sizes of rectangles of the Markov
partition $\cal R$. Similarly, if $C\in{\cal C}_2$ is defined on
an interval $\Lambda=[k^-,k^+]$ and $t=k_1^--k^->0$, then
$\pi(C)\subset R^s_{2,z}(\varepsilon)$ with
$\varepsilon=c\theta^t$.

Now define subcollections ${\cal C}_i'\subset {\cal C}_i$ for
$i=1,2$ that contain all cylinders $C$ defined on intervals
$\Lambda=[k^-,k^+]$ satisfying $k^+-k_1^+>k/3$ for $i=1$ and
$k_1^--k^->k/3$ for $i=2$ (recall that $k=\,{\rm
dist}(\Lambda_1,\Lambda_2)$). By the assumption (\ref{lne})
$$
   \mu\left (\cup_{C\in{\cal C}_i'} \pi(C)\right )
   \leq C_0\frac{\mu(R_i)}{|ak+b|^{1+\gamma}}
$$
with constants $a=-\ln\theta^{1/3}>0$ and $b=-\ln c$. So, the
parts $\pi(C)$, $C\in{\cal C}_i'$, can be removed from $R_i$ with
no harm. Denote by
$$
    \tilde{R}_i=R_i\setminus\left (\cup_{C\in{\cal C}_i'} \pi(C)\right )
$$
the remaining parts of $R_i$.

Note that $\tilde{R}_1$ and $\tilde{R}_2$ consist (mod 0) of
projections of cylinders $C'\in{\cal C}_1\setminus{\cal C}_1'$ and $C''\in{\cal
C}_2\setminus{\cal C}_2'$, respectively, and the gap between the
intervals on which $C'$ and $C''$ are defined is always $\geq
k/3$. Hence we can use the subadditivity of the correlation
function and Fact~3 of Section~\ref{secPTMC} to get
\begin{eqnarray*}
   |\mu(\tilde R_1\cap \tilde R_2)-\mu(\tilde R_1)\mu(\tilde R_2)|  &=&
   \big|\mu\big((\cup\pi(C')\big)\cap \big(\cup\pi(C'')\big)
   -\mu\big(\cup\pi(C')\big)\mu\big(\cup\pi(C'')\big)\big|\\
      &\leq &
   \sum_{C'}\sum_{C''}\big|\mu\big(\pi(C')\cap \pi(C'')\big)
   -\mu\big(\pi(C')\big)\mu(\pi(C'')\big)\big|
       \\ &\leq &
   \sum_{C'}\sum_{C''}
c_3\theta_3^{k/3}\mu\big(\pi(C')\big)\mu\big(\pi(C'')\big)
        \\ &\leq&
    c_3\theta_3^{k/3}\mu(\tilde{R}_1)\mu(\tilde{R}_2)\,.
\end{eqnarray*}
This completes the proof of Lemma~\ref{lmF3}. $\Box$

\begin{lemma}
Let $R_1,R_2$ be EQR rectangles with the corresponding intervals
$\Lambda_1$ and $\Lambda_2$. Then
$$
    |\mu(R_1\cap R_2)-\mu(R_1)\mu(R_2)|\leq
      c_9\frac{\mu(R_1)+\mu(R_2)}{|a_1\delta(\Lambda_1,\Lambda_2)+b|
^{1+\gamma}}\,.        
$$
Here $\delta(\Lambda_1,\Lambda_2)$ is the asymmetric distance defined
in Section 3, $a_1=a/2$, and $c_9>0$ is a constant.
\label{lmRR}
\end{lemma}

{\em Proof}. If $|\Lambda_2| \ge \frac12\delta(\Lambda_1,\Lambda_2)$, then,
by Lemmas~\ref{lmF1} and~\ref{lmF2}, both $\mu(R_1\cap R_2)$ and
$\mu(R_1)\mu(R_2)$ are bounded from above by $
c\theta^{\delta(\Lambda_1,\Lambda_2)}\mu(R_1)$, with  some $c>0$ and
$0 < \theta < 1$  depending only on $\mu$. Otherwise 
dist$(\Lambda_1,\Lambda_2)\geq \delta(\Lambda_1,\Lambda_2)/2$, and the
claim follows from 
Lemma~\ref{lmF3}. $\Box$\medskip

{\em Proof of Theorem~\ref{tm4}} \ goes by the same lines as the
proof of Theorem~\ref{tm1}. We estimate the quantity
$R_{mn}=\mu(R_m\cap T^{m-n}R_n)-\mu(R_n)\mu(R_m)$. Without loss of
generality, assume that  $\delta(\Lambda_n,\Lambda_m) \leq D$. By
Lemma~\ref{lmRR}
(applied to the rectangles $R_m$ and $T^{m-n}R_n$),
we have
$$
    |R_{mn}|\leq c_9\frac{\mu(R_m)+\mu(R_n)}{|a_1|m-n|-a_1D+b|^{1+\gamma}}\,.
$$
We use this bound if $|m-n|\geq D-b/a_1$, otherwise we can
use an obvious bound
$$
   |R_{mn}|\leq \mu(R_m)+\mu(R_n)\,.
$$
Summing up over all  over all $n$ with $\delta(\Lambda_n,\Lambda_m)
\leq D$, and then  over
$m=M,\ldots,N$, proves (SP). $\Box$\medskip

For the proof of Theorem~\ref{tm5}, we need two more lemmas. Recall
that now dim$\, X=2$ and every connected rectangle is EQR.

\begin{lemma}
Let $R_m$ and $R_{\ast}$ be two connected rectangles with $B$-bounded
aspect ratio, and $k^{\pm}_m$, $k^{\pm}_{\ast}$ integers characterizing
their quasiroundness. Assume that $R_m\cap R_{\ast}\neq\emptyset$.
If $d^u(R_m)\geq d^u(R_{\ast})$, then
$$
     k_{\ast}^+\geq k_m^++c_{10}\ln d^u(R_m)/d^u(R_{\ast})\,.
$$
Similarly, if $d^s(R_m)\geq d^s(R_{\ast})$, then
$$
     k_{\ast}^-\leq k_m^--c_{10}\ln d^s(R_m)/d^s(R_{\ast})\,.
$$
Here $c_{10}=c_{10}(\varepsilon_0,\varepsilon_1,B)>0$ is a constant.
\label{lmRRast}
\end{lemma}

{\em Proof}. This follows from standard distortion bounds. $\Box$ \medskip

{\em Proof of Theorem~\ref{tm5}}.
Denote by $k_n^{\pm}$ the integers characterizing
the quasiroundness of $R_n$. We may assume that all $R_n$
are small enough, and then the uniform boundedness of their
aspect ratio ensures that $k_n^+\geq 0$ and $k_n^-\leq 0$
for all $n\geq 1$.

We estimate the quantity
$R_{mn}=\mu(R_m\cap T^{m-n}R_n)-\mu(R_n)\mu(R_m)$.
The set $R_{\ast}=T^{m-n}R_n$ is a connected rectangle
whose quasiroundness is characterized by the integers
$k^{\pm}_{\ast}:=k^{\pm}_n+(n-m)$.
Without loss of generality, assume that $d^u(R_m)\geq d^u(R_n)$.

We consider three cases:

Case 1. Assume that either (i) $n-m\geq 2k_m^+$ or (ii)
$n-m\leq 2k_m^-$. In the case (i) we have
$$
      k^+_{\ast}-k^+_m\geq n-m-k_m^+\geq |n-m|/2\,,
$$
and in the case (ii) we have
$$
      k^-_m-k^-_{\ast}\geq k_m^--(n-m)\geq |n-m|/2\,.
$$
In either case we apply Lemma~\ref{lmRR} and obtain
$$
     |R_{mn}|\leq c_9\frac{\mu(R_m)+\mu(R_n)}{|a_2|n-m|+b|^{1+\gamma}}
$$
with $a_2=a_1/2>0$.

Case 2. Assume that $2k_m^-\leq n-m\leq 2k_m^+$ and $R_{\ast}\cap
R_m\neq\emptyset$.
If $n>m$, then $d^u(R_{\ast})\leq \theta^{n-m} d^u(R_m)$, and if
$n\leq m$, then $d^s(R_{\ast})\leq \theta^{m-n} d^s(R_m)$
for some constant $\theta<1$, due to the uniform hyperbolicity
of $T$. Hence, Lemma~\ref{lmRRast} implies that if $n>m$, then
$$
     k^+_{\ast}-k^+_m\geq c_{11}|n-m|\,
$$
and if $n\leq m$, then
$$
     k^-_m-k^-_{\ast}\geq c_{11}|n-m|
$$
with some constant $c_{11}>0$. Again, we use Lemma~\ref{lmRR} and obtain
$$
     |R_{mn}|\leq c_9\frac{\mu(R_m)+\mu(R_n)}{|a_3|n-m|+b|^{1+\gamma}}
$$
with $a_3=c_{11}a_1$.

Case 3. Assume that $2k_m^-\leq n-m\leq 2k_m^+$ and $R_{\ast}\cap
R_m=\emptyset$. Then $R_{mn}=\mu(R_m)\mu(R_n)$.
It follows from Lemma~\ref{lmF1} that
$$
       a_4|\ln\mu(R_m)|+b_4\leq k^+_m-k^-_m\leq a_5|\ln\mu(R_m)|+b_5
$$
with some $a_4,a_5>0$ and $-\infty<b_4,b_5<\infty$, and similar
bounds hold for $R_n$. Our assumption $d^u(R_n)\leq d^u(R_m)$
and the $B$-boundedness of aspect ratio imply that
$k_n^+\geq \varepsilon_2 k_m^+$ and
$k_n^-\leq \varepsilon_2 k_m^-$ for some constant $\varepsilon_2>0$,
due to uniform bounds on expansion and contraction rates of $T$.
Therefore,
$$
      \mu(R_n)\leq c_{12}[\mu(R_m)]^{\kappa}
$$
with some constants $c_{12}>0$ and $\kappa>0$. Holding
$m$ fixed and summing over all $n$ that satisfy the conditions
of Case 3 gives
$$
  \sum_n R_{mn}\leq 2c_{12}[\mu(R_m)]^{1+\kappa}(a_5|\ln\mu(R_m)|+b_5)
     \leq c_{13}\mu(R_m)
$$
with some constant $c_{13}>0$.

Lastly, summing up over all $m,n=M,\ldots,N$ proves (SP).
$\Box$\bigskip

\bigskip

{\bf Acknowledgements}. The authors want to thank Vitaly
Bergelson, Jean-Pierre Conze,  Dmitry Dolgopyat, Yves Guivarc'h,
Andres del Junco and Albert Raugi for helpful discussions, and the
referee for useful comments.


\begin{thebibliography}{99}

\bibitem{Bos} M. Boshernitzan, {\em Quantitative recurrence
results}, Invent. Math. {\bf 113} (1993), 617--631.

\bibitem{Bo75} R. Bowen, {\em Equilibrium states and
the ergodic theory of Anosov diffeomorphisms},
Lect. Notes Math. {\bf 470}, Springer-Verlag,
Berlin, 1975.

\bibitem{ChS} N. Chernov, {\em Invariant measures
for hyperbolic dynamical systems}, to appear in In: {\em Handbook
of Dynamical Systems, Vol. I}, Ed. A.~Katok and B.~Hasselblatt,
Elsevier.


\bibitem{CR}  J.-P. Conze and A. Raugi, {\em Convergence des
potentiels pour un op\'erateur de transfert,
applications aux syst\`emes dynamiques et aux cha\^\i nes de
Markov}, manuscript.

\bibitem{Do} D. Dolgopyat, {\em Limit theorems for partially
hyperbolic systems}, manuscript.

\bibitem{Durr} R. Durrett, {\em Probability: theory and examples},
Wadsworth \& Brooks/Cole, 1991

\bibitem{FT} N. Friedman and E. Thomas, {\em Higher order sweeping
out}, Ill. J. Math. {\bf 29} (1985), 401--417.

\bibitem{HV} R. Hill and S. Velani, {\em  Ergodic theory of
shrinking targets},    Invent. Math. {\bf 119} (1995),
175--198.

\bibitem{KM} D. Kleinbock and G. Margulis, {\em Logarithm laws for
flows on homogeneous spaces}, Inv.~Math.~{\bf 138} (1999), 451--494.

 \bibitem{Kr} U. Krengel, {\em On the individual ergodic theorem for
subsequences}, Ann. of
Math. Stat. {\bf 42} (1971), 1091-1095.

\bibitem{Ph} W. Philipp, {\em Some metrical theorems in number
theory}, Pacific J. Math. {\bf 20} (1967), 109--127.

\bibitem{Ru99} D. Ruelle, {\em Smooth dynamics and new theoretical
ideas in nonequilibrium statistical mechanics}, J. Statist. Phys.
{\bf 95} (1999), 393--468.

\bibitem{Pet} K. Petersen, {\em Ergodic theory}, Camb. Univ.
Press, 1983.

\bibitem{Sp} V. Sprind\v{z}uk, {\em Metric theory of Diophantine
approximations}, J. Wiley \& Sons, New York-Toronto-London, 1979.

\end{thebibliography}
\end{document}